\documentclass[11pt,reqno]{amsart}

\usepackage{amssymb}
\usepackage{amscd}
\usepackage{amsfonts}
\usepackage{setspace}
\usepackage{version}
\usepackage{hyperref}

\usepackage{graphicx}

\theoremstyle{definition}

\renewcommand{\leq}{\leqslant}
\renewcommand{\geq}{\geqslant}

\def\Z{\mathbf{Z}}

\newcommand{\md}[1]{\ensuremath{(\operatorname{mod}\, #1)}}

\numberwithin{equation}{section}

\begin{document}

\title[Lower bounds for corner-free sets]{Lower bounds for corner-free sets}


\author{Ben Green}
\address{Mathematical Institute\\
Radcliffe Observatory Quarter\\
Woodstock Road\\
Oxford OX2 6GG\\
England }
\email{ben.green@maths.ox.ac.uk}


\thanks{The author is supported by a Simons Investigator Grant.}

\begin{abstract}
A \emph{corner} is a set of three points in $\Z^2$ of the form $(x, y), (x + d, y), (x, y + d)$ with $d \neq 0$. We show that for infinitely many $N$ there is a set $A \subset [N]^2$ of size $2^{-(c + o(1)) \sqrt{\log_2 N}} N^2$ not containing any corner, where $c = 2 \sqrt{2 \log_2 \frac{4}{3}} \approx 1.822\dots$.
\end{abstract}
\maketitle

Let $q,d$ be large positive integers. For each $x \in [q^d -1]$, write $\pi(x) = (x_0,\dots, x_{d-1}) \in \Z^d$ for the vector of digits of its base $q$ expansion, thus $x = \sum_{i=0}^{d-1} x_i q^{i}$, with $0 \leq x_i < q$ for all $i$. 

For each positive integer $r$, consider the set $A_r$ of all pairs $(x,y) \in [q^d -1]^2$ for which $\Vert \pi(x) - \pi(y) \Vert_2^2 = r$ and $\frac{q}{2} \leq x_i + y_i < \frac{3q}{2}$ for all $i$.

We claim that $A_r$ is free of corners.  Suppose that $(x,y), (x+d, y), (x, y + d) \in A_r$.  Then
\begin{equation}\label{sphere} \Vert \pi(x) - \pi(y) \Vert^2_2 = \Vert \pi(x+d) - \pi(y)\Vert^2_2 = \Vert \pi(x) - \pi(y + d) \Vert^2_2 = r.\end{equation}
We claim that 
\begin{equation}\label{cocycle} \pi(x + d) + \pi(y) = \pi(x) + \pi(y + d).\end{equation} To this end, we show that $(x+d)_i + y_i = x_i + (y+d)_i$ for $i = 0,1,\dots$ by induction on $i$. A single argument works for both the base case $i = 0$ and the inductive step. Suppose that, for some $j \geq 0$, we have the statement for $i < j$. Write $x_{\geq j} := \sum_{i \geq  j} x_i q^i$, and define $(x + d)_{\geq j}, y_{\geq j}, (y+d)_{\geq j}$ similarly. By the inductive hypothesis and the fact that $x + (y + d) = (x + d) + y$, we see that $x_{\geq j} + (y + d)_{\geq j} = (x + d)_{\geq j} + y_{\geq j}$. Therefore $x_{ j} + (y + d)_{ j} = (x + d)_{ j} + y_{ j} \md{q}$. However by assumption we have $\frac{q}{2} \leq x_{ j} + (y + d)_{ j} , (x + d)_{ j} + y_{ j} < \frac{3q}{2}$, and so $x_{ j} + (y + d)_{ j} = (x + d)_{ j} + y_{ j}$. The induction goes through.

With \eqref{cocycle} established, let us return to \eqref{sphere}. We now see that this statement implies that $\Vert a  \Vert_2^2 = \Vert a + b \Vert_2^2 = \Vert a - b \Vert_2^2 = r$, where $a := \pi(x) - \pi(y)$ and $b := \pi(x + d) - \pi(x) = \pi(y + d) - \pi(y)$. By the parallelogram law $2 \Vert a \Vert_2^2 + 2 \Vert b \Vert_2^2 = \Vert a - b \Vert_2^2 + \Vert a + b \Vert_2^2$, this immediately implies that $b = 0$. Since $\pi$ is injective, it follows that $d = 0$ and so indeed $A_r$ is corner-free.

The set of all pairs $(x,y)$ with $\frac{q}{2} \leq x_i + y_i < \frac{3q}{2}$ for all $i$ has size $(\frac{3}{4}q^2 + O(q))^d$. Therefore by the pigeonhole principle there is some $r$ such that $\# A_r \geq (dq^2)^{-1} (\frac{3}{4}q^2 + O(q))^d$. 

Now for a given $d$ set $q := \lfloor (2/\sqrt{3})^d \rfloor$ and $N := q^d$. Then $A_r \subset [N] \times [N]$, $A_r$ is free of corners, and
\[ \# A_r  \geq N^2 (dq^2)^{-1}(\frac{3}{4} + O(\frac{1}{q}))^d .\]
Writing $o(1)$ for a quantity tending to $0$ as $N \rightarrow \infty$, we note that $q = (\frac{2}{\sqrt{3}} + o(1))^d$ and that $d = (1 + o(1))\sqrt{\frac{\log_2 N}{\log_2(2/\sqrt{3})}}$. A short calculation then confirms that 
\[ \# A_r \geq  N^2 2^{-(c + o(1)) \sqrt{\log_2 N}},\]
where $c = 2 \sqrt{2 \log_2 \frac{4}{3}} \approx 1.822\dots$. \vspace*{7pt}

\emph{Remark.} The construction came about by a careful study of the recent preprint of Linial and Shraibman \cite{linial-shraibman}, where they used ideas from communication complexity to obtain a bound with $c = 2\sqrt{\log_2 e} \approx 2.402
\dots$, improving on the previously best known bound with $c = 2 \sqrt{2} \approx 2.828\dots$ which comes from Behrend's construction. By bypassing the language of communication complexity one may simplify the construction, in particular avoiding the use of entropy methods. This yields a superior bound.


\begin{thebibliography}{99}
\bibitem{linial-shraibman} N.~Linial and A.~Shraibman, \emph{Larger corner-free sets from better NOF exactly-$N$ protocols,} preprint, arxiv:2102.00421.
\end{thebibliography}
\end{document}